# Optimal sampling strategies for multiscale stochastic processes

Vinay J. Ribeiro[1], Rudolf H. Riedi[1] and Richard G. Baraniuk[2],*

*Rice University*

**Abstract:** In this paper, we determine which non-random sampling of fixed size gives the best linear predictor of the sum of a finite spatial population. We employ different multiscale superpopulation models and use the minimum mean-squared error as our optimality criterion. In multiscale superpopulation tree models, the leaves represent the units of the population, interior nodes represent partial sums of the population, and the root node represents the total sum of the population. We prove that the optimal sampling pattern varies dramatically with the correlation structure of the tree nodes. While *uniform sampling* is optimal for trees with "positive correlation progression", it provides the worst possible sampling with "negative correlation progression." As an analysis tool, we introduce and study a class of *independent innovations trees* that are of interest in their own right. We derive a fast water-filling algorithm to determine the optimal sampling of the leaves to estimate the root of an independent innovations tree.

## 1. Introduction

In this paper we design optimal sampling strategies for spatial populations under different multiscale superpopulation models. Spatial sampling plays an important role in a number of disciplines, including geology, ecology, and environmental science. See, e.g., Cressie [5].

### 1.1. Optimal spatial sampling

Consider a finite population consisting of a rectangular grid of $R \times C$ units as depicted in Fig. 1(a). Associated with the unit in the $i^{\text{th}}$ row and $j^{\text{th}}$ column is an unknown value $\ell_{i,j}$. We treat the $\ell_{i,j}$'s as one realization of a superpopulation model.

Our goal is to determine which sample, among all samples of size $n$, gives the best linear estimator of the population sum, $S := \sum_{i,j} \ell_{i,j}$. We abbreviate *variance*, *covariance*, and *expectation* by "var", "cov", and "$\mathbb{E}$" respectively. Without loss of generality we assume that $\mathbb{E}(\ell_{i,j}) = 0$ for all locations $(i, j)$.

[1]Department of Statistics, 6100 Main Street, MS-138, Rice University, Houston, TX 77005, e-mail: vinay@rice.edu; riedi@rice.edu
[2]Department of Electrical and Computer Engineering, 6100 Main Street, MS-380, Rice University, Houston, TX 77005, e-mail: richb@rice.edu, url: dsp.rice.edu, spin.rice.edu
*Supported by NSF Grants ANI-9979465, ANI-0099148, and ANI-0338856, DoE SciDAC Grant DE-FC02-01ER25462, DARPA/AFRL Grant F30602-00-2-0557, Texas ATP Grant 003604-0036-2003, and the Texas Instruments Leadership University program.

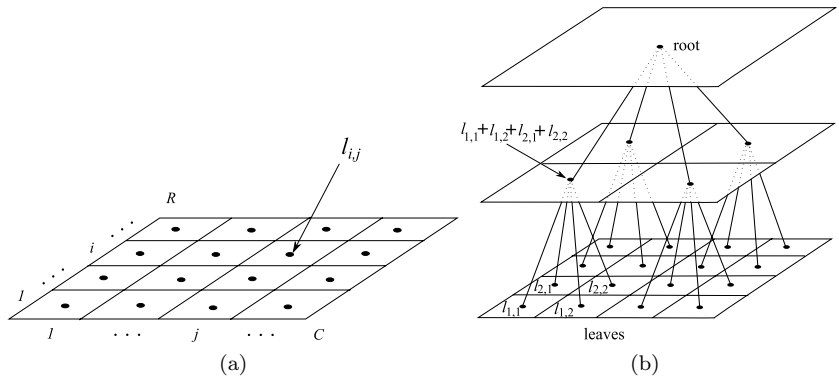

FIG 1. (a) *Finite population on a spatial rectangular grid of size $R \times C$ units. Associated with the unit at position $(i,j)$ is an unknown value $\ell_{i,j}$.* (b) *Multiscale superpopulation model for a finite population. Nodes at the bottom are called leaves and the topmost node the root. Each leaf node corresponds to one value $\ell_{i,j}$. All nodes, except for the leaves, correspond to the sum of their children at the next lower level.*

Denote an arbitrary sample of size $n$ by $L$. We consider linear estimators of $S$ that take the form

$$\widehat{S}(L, \boldsymbol{\alpha}) := \boldsymbol{\alpha}^T L, \tag{1.1}$$

where $\boldsymbol{\alpha}$ is an arbitrary set of coefficients. We measure the accuracy of $\widehat{S}(L, \boldsymbol{\alpha})$ in terms of the *mean-squared error* (MSE)

$$\mathcal{E}(S|L, \boldsymbol{\alpha}) := \mathbb{E}\left(S - \widehat{S}(L, \boldsymbol{\alpha})\right)^2 \tag{1.2}$$

and define the *linear minimum mean-squared error* (LMMSE) of estimating $S$ from $L$ as

$$\mathcal{E}(S|L) := \min_{\boldsymbol{\alpha} \in \mathbb{R}^n} \mathcal{E}(S|L, \boldsymbol{\alpha}). \tag{1.3}$$

Restated, our goal is to determine

$$L^* := \arg\min_L \mathcal{E}(S|L). \tag{1.4}$$

Our results are particularly applicable to Gaussian processes for which linear estimation is optimal in terms of mean-squared error. We note that for certain multimodal and discrete processes linear estimation may be sub-optimal.

### 1.2. Multiscale superpopulation models

We assume that the population is one realization of a multiscale stochastic process (see Fig. 1(b)) (see Willsky [20]). Such processes consist of random variables organized on a tree. Nodes at the bottom, called *leaves*, correspond to the population $\ell_{i,j}$. All nodes, except for the leaves, represent the sum total of their children at the next lower level. The topmost node, the *root*, hence represents the sum of the entire population. The problem we address in this paper is thus equivalent to the following: *Among all possible sets of leaves of size $n$, which set provides the best linear estimator for the root in terms of MSE?*

Multiscale stochastic processes efficiently capture the correlation structure of a wide range of phenomena, from uncorrelated data to complex *fractal* data. They do so through a simple probabilistic relationship between each parent node and its children. They also provide fast algorithms for analysis and synthesis of data and are often physically motivated. As a result multiscale processes have been used in a number of fields, including oceanography, hydrology, imaging, physics, computer networks, and sensor networks (see Willsky [20] and references therein, Riedi et al. [15], and Willett et al. [19]).

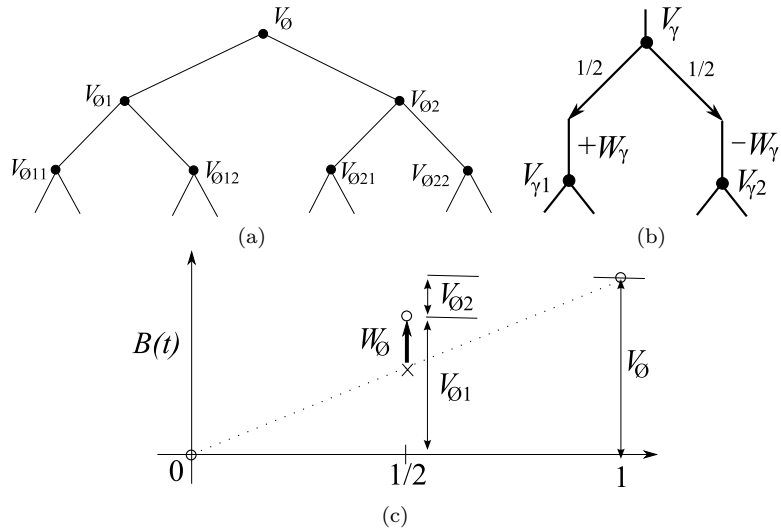

FIG 2. (a) *Binary tree for interpolation of Brownian motion, $B(t)$.* (b) *Form child nodes $V_{\gamma 1}$ and $V_{\gamma 2}$ by adding and subtracting an independent Gaussian random variable $W_\gamma$ from $V_\gamma/2$.* (c) *Mid-point displacement. Set $B(1) = V_\emptyset$ and form $B(1/2) = (B(1) - B(0))/2 + W_\emptyset = V_{\emptyset 1}$. Then $B(1) - B(1/2) = V_\emptyset/2 - W_\emptyset = V_{\emptyset 2}$. In general a node at scale $j$ and position $k$ from the left of the tree corresponds to $B((k+1)2^{-j}) - B(k2^{-j})$.*

We illustrate the essentials of multiscale modeling through a tree-based interpolation of one-dimensional *standard Brownian motion*. Brownian motion, $B(t)$, is a zero-mean Gaussian process with $B(0) := 0$ and $\text{var}(B(t)) = t$. Our goal is to begin with $B(t)$ specified only at $t = 1$ and then interpolate it at all time instants $t = k2^{-j}$, $k = 1, 2, \ldots, 2^j$ for any given value $j$.

Consider a binary tree as shown in Fig. 2(a). We denote the root by $V_\emptyset$. Each node $V_\gamma$ is the *parent* of two nodes connected to it at the next lower level, $V_{\gamma 1}$ and $V_{\gamma 2}$, which are called its *child* nodes. The address $\gamma$ of any node $V_\gamma$ is thus a concatenation of the form $\emptyset k_1 k_2 \ldots k_j$, where $j$ is the node's *scale* or depth in the tree.

We begin by generating a zero-mean Gaussian random variable with unit variance and assign this value to the root, $V_\emptyset$. The root is now a realization of $B(1)$. We next interpolate $B(0)$ and $B(1)$ to obtain $B(1/2)$ using a "mid-point displacement" technique. We generate independent *innovation* $W_\emptyset$ of variance $\text{var}(W_\emptyset) = 1/4$ and set $B(1/2) = V_\emptyset/2 + W_\emptyset$ (see Fig. 2(c)).

Random variables of the form $B((k+1)2^{-j}) - B(k2^{-j})$ are called *increments* of Brownian motion at time-scale $j$. We assign the increments of the Brownian motion at time-scale 1 to the children of $V_\emptyset$. That is, we set

$$
\begin{aligned}
V_{\emptyset 1} &= B(1/2) - B(0) = V_\emptyset/2 + W_\emptyset, \text{ and} \\
V_{\emptyset 2} &= B(1) - B(1/2) = V_\emptyset/2 - W_\emptyset
\end{aligned}
\tag{1.5}
$$

as depicted in Fig. 2(c). We continue the interpolation by repeating the procedure described above, replacing $V_\emptyset$ by each of its children and reducing the variance of the innovations by half, to obtain $V_{\emptyset 11}$, $V_{\emptyset 12}$, $V_{\emptyset 21}$, and $V_{\emptyset 22}$.

Proceeding in this fashion we go down the tree assigning values to the different tree nodes (see Fig. 2(b)). It is easily shown that the nodes at scale $j$ are now realizations of $B((k+1)2^{-j}) - B(k2^{-j})$. That is, increments at time-scale $j$. For a given value of $j$ we thus obtain the interpolated values of Brownian motion, $B(k2^{-j})$ for $k = 0, 1, \ldots, 2^j - 1$, by cumulatively summing up the nodes at scale $j$.

By appropriately setting the variances of the innovations $W_\gamma$, we can use the procedure outlined above for Brownian motion interpolation to interpolate several other Gaussian processes (Abry et al. [1], Ma and Ji [12]). One of these is *fractional Brownian motion* (fBm), $B_H(t)$ ($0 < H < 1$), that has variance $\text{var}(B_H(t)) = t^{2H}$. The parameter $H$ is called the *Hurst* parameter. Unlike the interpolation for Brownian motion which is exact, however, the interpolation for fBm is only approximate. By setting the variance of innovations at different scales appropriately we ensure that nodes at scale $j$ have the same variance as the increments of fBm at time-scale $j$. However, except for the special case when $H = 1/2$, the covariance between any two arbitrary nodes at scale $j$ is not always identical to the covariance of the corresponding increments of fBm at time-scale $j$. Thus the tree-based interpolation captures the variance of the increments of fBm at all time-scales $j$ but does not perfectly capture the entire covariance (second-order) structure.

This approximate interpolation of fBm, nevertheless, suffices for several applications including network traffic synthesis and queuing experiments (Ma and Ji [12]). They provide fast $O(N)$ algorithms for both synthesis and analysis of data sets of size $N$. By assigning multivariate random variables to the tree nodes $V_\gamma$ as well as innovations $W_\gamma$, the accuracy of the interpolations for fBm can be further improved (Willsky [20]).

In this paper we restrict our attention to two types of multiscale stochastic processes: *covariance trees* (Ma and Ji [12], Riedi et al. [15]) and *independent innovations trees* (Chou et al. [3], Willsky [20]). In covariance trees the covariance between pairs of leaves is purely a function of their distance. In independent innovations trees, each node is related to its parent nodes through a unique independent additive innovation. One example of a covariance tree is the multiscale process described above for the interpolation of Brownian motion (see Fig. 2).

### *1.3. Summary of results and paper organization*

We analyze covariance trees belonging to two broad classes: those with *positive correlation progression* and those with *negative correlation progression*. In trees with positive correlation progression, leaves closer together are more correlated than leaves father apart. The opposite is true for trees with negative correlation progression. While most spatial data sets are better modeled by trees with positive correlation progression, there exist several phenomena in finance, computer networks, and nature that exhibit anti-persistent behavior, which is better modeled by a tree with negative correlation progression (Li and Mills [11], Kuchment and Gelfan [9], Jamdee and Los [8]).

For covariance trees with positive correlation progression we prove that uniformly spaced leaves are optimal and that clustered leaf nodes provides the worst possible MSE among all samples of fixed size. The optimal solution can, however, change with the correlation structure of the tree. In fact for covariance trees with negative correlation progression we prove that uniformly spaced leaf nodes give the *worst* possible MSE!

In order to prove optimality results for covariance trees we investigate the closely related independent innovations trees. In these trees, a parent node cannot equal the sum of its children. As a result they cannot be used as superpopulation models in the scenario described in Section 1.1. Independent innovations trees are however of interest in their own right. For independent innovations trees we describe an efficient algorithm to determine an optimal leaf set of size $n$ called *water-filling*. Note that the general problem of determining which $n$ random variables from a given set provide the best linear estimate of another random variable that is not in the same set is an NP-hard problem. In contrast, the water-filling algorithm solves one problem of this type in polynomial-time.

The paper is organized as follows. Section 2 describes various multiscale stochastic processes used in the paper. In Section 3 we describe the water-filling technique to obtain optimal solutions for independent innovations trees. We then prove optimal and worst case solutions for covariance trees in Section 4. Through numerical experiments in Section 5 we demonstrate that optimal solutions for multiscale processes can vary depending on their topology and correlation structure. We describe related work on optimal sampling in Section 6. We summarize the paper and discuss future work in Section 7. The proofs can be found in the Appendix. The pseudo-code and analysis of the computational complexity of the water-filling algorithm are available online (Ribeiro et al. [14]).

## 2. Multiscale stochastic processes

Trees occur naturally in many applications as an efficient data structure with a simple dependence structure. Of particular interest are trees which arise from representing and analyzing stochastic processes and time series on different time scales. In this section we describe various trees and related background material relevant to this paper.

### 2.1. Terminology and notation

A tree is a special graph, i.e., a set of nodes together with a list of pairs of nodes which can be pictured as directed edges pointing from one node to another with the following special properties (see Fig. 3): (1) There is a unique node called the *root* to which no edge points to. (2) There is exactly one edge pointing to any node, with the exception of the root. The starting node of the edge is called the *parent* of the ending node. The ending node is called a *child* of its parent. (3) The tree is *connected*, meaning that it is possible to reach any node from the root by following edges.

These simple rules imply that there are no cycles in the tree, in particular, there is exactly one way to reach a node from the root. Consequently, unique addresses can be assigned to the nodes which also reflect the level of a node in the tree. The topmost node is the root whose address we denote by ø. Given an arbitrary node $\gamma$, its child nodes are said to be one level lower in the tree and are addressed by $\gamma k$ ($k = 1, 2, \ldots, P_\gamma$), where $P_\gamma \geq 0$. The address of each node is thus a concatenation of the form $ø k_1 k_2 \ldots k_j$, or $k_1 k_2 \ldots k_j$ for short, where $j$ is the node's *scale* or depth in the tree. The largest scale of any node in the tree is called the *depth* of the tree.

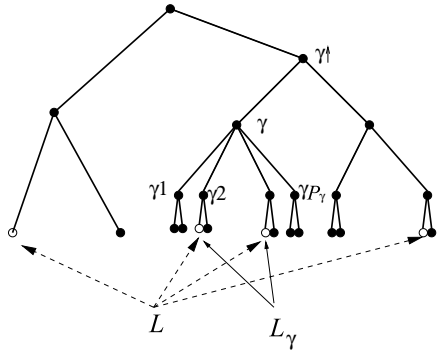

Fig 3. *Notation for multiscale stochastic processes.*

Nodes with no child nodes are termed *leaves* or *leaf nodes*. As usual, we denote the number of elements of a set of leaf nodes $L$ by $|L|$. We define the operator $\uparrow$ such that $\gamma k \uparrow = \gamma$. Thus, the operator $\uparrow$ takes us one level higher in the tree to the parent of the current node. Nodes that can be reached from $\gamma$ by repeated $\uparrow$ operations are called *ancestors* of $\gamma$. We term $\gamma$ a *descendant* of all of its ancestors.

The set of nodes and edges formed by $\gamma$ and all its descendants is termed the *tree of* $\gamma$. Clearly, it satisfies all rules of a tree. Let $L_\gamma$ denote the subset of $L$ that belong to the tree of $\gamma$. Let $\mathcal{N}_\gamma$ be the total number of leaves of the tree of $\gamma$.

To every node $\gamma$ we associate a single (univariate) random variable $V_\gamma$. For the sake of brevity we often refer to $V_\gamma$ as simply "the node $V_\gamma$" rather than "the random variable associated with node $\gamma$."

## 2.2. Covariance trees

Covariance trees are multiscale stochastic processes defined on the basis of the covariance between the leaf nodes which is purely a function of their *proximity*. Examples of covariance trees are the Wavelet-domain Independent Gaussian model (WIG) and the Multifractal Wavelet Model (MWM) proposed for network traffic (Ma and Ji [12], Riedi et al. [15]). Precise definitions follow.

**Definition 2.1.** The proximity of two leaf nodes is the scale of their lowest common ancestor.

Note that the larger the proximity of a pair of leaf nodes, the closer the nodes are to each other in the tree.

**Definition 2.2.** A covariance tree is a multiscale stochastic process with two properties. (1) The covariance of any two leaf nodes depends only on their proximity. In other words, if the leaves $\gamma'$ and $\gamma$ have proximity $k$ then $\text{cov}(V_\gamma, V_{\gamma'}) =: c_k$. (2) All leaf nodes are at the same scale $D$ and the root is equally correlated with all leaves.

In this paper we consider covariance trees of two classes: trees with positive correlation progression and trees with negative correlation progression.

**Definition 2.3.** A covariance tree has a *positive correlation progression* if $c_k > c_{k-1} > 0$ for $k = 1, \ldots, D-1$. A covariance tree has a *negative correlation progression* if $c_k < c_{k-1}$ for $k = 1, \ldots, D-1$.

Intuitively in trees with positive correlation progression leaf nodes "closer" to each other in the tree are more strongly correlated than leaf nodes "farther apart." Our results take on a special form for covariance trees that are also symmetric trees.

**Definition 2.4.** A symmetric tree is a multiscale stochastic process in which $P_\gamma$, the number of child nodes of $V_\gamma$, is purely a function of the scale of $\gamma$.

### 2.3. Independent innovations trees

Independent innovations trees are particular multiscale stochastic processes defined as follows.

**Definition 2.5.** An independent innovations tree is a multiscale stochastic process in which each node $V_\gamma$, excluding the root, is defined through

$$(2.1) \qquad V_\gamma := \varrho_\gamma V_{\gamma\uparrow} + W_\gamma.$$

Here, $\varrho_\gamma$ is a scalar and $W_\gamma$ is a random variable independent of $V_{\gamma\uparrow}$ as well as of $W_{\gamma'}$ for all $\gamma' \neq \gamma$. The root, $V_\emptyset$, is independent of $W_\gamma$ for all $\gamma$. In addition $\varrho_\gamma \neq 0$, $\text{var}(W_\gamma) > 0 \ \forall \gamma$ and $\text{var}(V_\emptyset) > 0$.

Note that the above definition guarantees that $\text{var}(V_\gamma) > 0 \ \forall \gamma$ as well as the *linear independence*[1] of any set of tree nodes.

The fact that each node is the sum of a scaled version of its parent and an independent random variable makes these trees amenable to analysis (Chou et al. [3], Willsky [20]). We prove optimality results for independent innovations trees in Section 3. Our results take on a special form for scale-invariant trees defined below.

**Definition 2.6.** A scale-invariant tree is an independent innovations tree which is symmetric and where $\varrho_\gamma$ and the distribution of $W_\gamma$ are purely functions of the scale of $\gamma$.

While independent innovations trees are not covariance trees in general, it is easy to see that scale-invariant trees are indeed covariance trees with positive correlation progression.

### 3. Optimal leaf sets for independent innovations trees

In this section we determine the optimal leaf sets of independent innovations trees to estimate the root. We first describe the concept of water-filling which we later use to prove optimality results. We also outline an efficient numerical method to obtain the optimal solutions.

#### 3.1. Water-filling

While obtaining optimal sets of leaves to estimate the root we maximize a sum of concave functions under certain constraints. We now develop the tools to solve this problem.

---

[1] A set of random variables is linearly independent if none of them can be written as a linear combination of finitely many other random variables in the set.

**Definition 3.1.** A real function $\psi$ defined on the set of integers $\{0, 1, \ldots, M\}$ is discrete-concave if

(3.1) $\quad\quad \psi(x+1) - \psi(x) \geq \psi(x+2) - \psi(x+1),\ \text{for}\ x = 0, 1, \ldots, M-2.$

The optimization problem we are faced with can be cast as follows. Given integers $P \geq 2$, $M_k > 0$ ($k = 1, \ldots, P$) and $n \leq \sum_{k=1}^{P} M_k$ consider the discrete space

(3.2) $\quad \Delta_n(M_1, \ldots, M_P) := \left\{ X = [x_k]_{k=1}^{P} : \sum_{k=1}^{P} x_k = n; x_k \in \{0, 1, \ldots, M_k\}, \forall k \right\}.$

Given non-decreasing, discrete-concave functions $\psi_k$ ($k = 1, \ldots, P$) with domains $\{0, \ldots, M_k\}$ we are interested in

(3.3) $\quad\quad h(n) := \max \left\{ \sum_{k=1}^{P} \psi_k(x_k)\ :\ X \in \Delta_n(M_1, \ldots, M_P) \right\}.$

In the context of optimal estimation on a tree, $P$ will play the role of the number of children that a parent node $V_\gamma$ has, $M_k$ the total number of leaf node descendants of the $k$-th child $V_{\gamma k}$, and $\psi_k$ the reciprocal of the optimal LMMSE of estimating $V_\gamma$ given $x_k$ leaf nodes in the tree of $V_{\gamma k}$. The quantity $h(n)$ corresponds to the reciprocal of the optimal LMMSE of estimating node $V_\gamma$ given $n$ leaf nodes in its tree.

The following iterative procedure solves the optimization problem (3.3). Form vectors $G^{(n)} = [g_k^{(n)}]_{k=1}^{P}$, $n = 0, \ldots, \sum_k M_k$ as follows:
Step (i): Set $g_k^{(0)} = 0$, $\forall k$.
Step (ii): Set

(3.4) $\quad\quad\quad\quad\quad\quad g_k^{(n+1)} = \begin{cases} g_k^{(n)} + 1, & k = m \\ g_k^{(n)}, & k \neq m \end{cases}$

where

(3.5) $\quad\quad m \in \arg\max_{k} \left\{ \psi_k\left(g_k^{(n)} + 1\right) - \psi_k\left(g_k^{(n)}\right) : g_k^{(n)} < M_k \right\}.$

The procedure described in Steps (i) and (ii) is termed *water-filling* because it resembles the solution to the problem of filling buckets with water to maximize the sum of the heights of the water levels. These buckets are narrow at the bottom and monotonically widen towards the top. Initially all buckets are empty (compare Step (i)). At each step we are allowed to pour one unit of water into any one bucket with the goal of maximizing the sum of water levels. Intuitively at any step we must pour the water into that bucket which will give the maximum increase in water level among all the buckets not yet full (compare Step (ii)). Variants of this water-filling procedure appear as solutions to different information theoretic and communication problems (Cover and Thomas [4]).

**Lemma 3.1.** *The function $h(n)$ is non-decreasing and discrete-concave. In addition,*

(3.6) $\quad\quad\quad\quad\quad\quad h(n) = \sum_{k} \psi_k\left(g_k^{(n)}\right),$

*where $g_k^{(n)}$ is defined through water-filling.*

When all functions $\psi_k$ in Lemma 3.1 are identical, the maximum of $\sum_{k=1}^{P} \psi_k(x_k)$ is achieved by choosing the $x_k$'s to be "near-equal". The following Corollary states this rigorously.

**Corollary 3.1.** *If $\psi_k = \psi$ for all $k = 1, 2, \ldots, P$ with $\psi$ non-decreasing and discrete-concave, then*

$$(3.7) \qquad h(n) = \left(P - n + P\left\lfloor\frac{n}{P}\right\rfloor\right)\psi\left(\left\lfloor\frac{n}{P}\right\rfloor\right) + \left(n - P\left\lfloor\frac{n}{P}\right\rfloor\right)\psi\left(\left\lfloor\frac{n}{P}\right\rfloor + 1\right).$$

*The maximizing values of the $x_k$ are apparent from (3.7). In particular, if $n$ is a multiple of $P$ then this reduces to*

$$(3.8) \qquad h(n) = P\psi\left(\frac{n}{P}\right).$$

Corollary 3.1 is key to proving our results for scale-invariant trees.

### 3.2. Optimal leaf sets through recursive water-filling

Our goal is to determine a choice of $n$ leaf nodes that gives the smallest possible LMMSE of the root. Recall that the LMMSE of $V_\gamma$ given $L_\gamma$ is defined as

$$(3.9) \qquad \mathcal{E}(V_\gamma | L_\gamma) := \min_{\boldsymbol{\alpha}} \mathbb{E}(V_\gamma - \boldsymbol{\alpha}^T L_\gamma)^2,$$

where, in an abuse of notation, $\boldsymbol{\alpha}^T L_\gamma$ denotes a linear combination of the elements of $L_\gamma$ with coefficients $\boldsymbol{\alpha}$. Crucial to our proofs is the fact that (Chou et al. [3] and Willsky [20]),

$$(3.10) \qquad \frac{1}{\mathcal{E}(V_\gamma | L_\gamma)} + \frac{P_\gamma - 1}{\operatorname{var}(V_\gamma)} = \sum_{k=1}^{P_\gamma} \frac{1}{\mathcal{E}(V_\gamma | L_{\gamma k})}.$$

Denote the set consisting of all subsets of leaves of the tree of $\gamma$ of size $n$ by $\Lambda_\gamma(n)$. Motivated by (3.10) we introduce

$$(3.11) \qquad \mu_\gamma(n) := \max_{L \in \Lambda_\gamma(n)} \mathcal{E}(V_\gamma | L)^{-1}$$

and define

$$(3.12) \qquad \mathcal{L}_\gamma(n) := \{L \in \Lambda_\gamma(n) : \mathcal{E}(V_\gamma | L)^{-1} = \mu_\gamma(n)\}.$$

Restated, our goal is to determine one element of $\mathcal{L}_\emptyset(n)$. To allow a recursive approach through scale we generalize (3.11) and (3.12) by defining

$$(3.13) \qquad \mu_{\gamma,\gamma'}(n) := \max_{L \in \Lambda_{\gamma'}(n)} \mathcal{E}(V_\gamma | L)^{-1} \quad \text{and}$$

$$(3.14) \qquad \mathcal{L}_{\gamma,\gamma'}(n) := \{L \in \Lambda_{\gamma'}(n) : \mathcal{E}(V_\gamma | L)^{-1} = \mu_{\gamma,\gamma'}(n)\}.$$

Of course, $\mathcal{L}_\gamma(n) = \mathcal{L}_{\gamma,\gamma}(n)$. For the recursion, we are mostly interested in $\mathcal{L}_{\gamma,\gamma k}(n)$, i.e., the optimal estimation of a parent node from a sample of leaf nodes of one of its children. The following will be useful notation

$$(3.15) \qquad X^* = [x_k^*]_{k=1}^{P_\gamma} := \arg\max_{X \in \Delta_n(\mathcal{N}_{\gamma 1}, \ldots, \mathcal{N}_{\gamma P_\gamma})} \sum_{k=1}^{P_\gamma} \mu_{\gamma,\gamma k}(x_k).$$

Using (3.10) we can decompose the problem of determining $L \in \mathcal{L}_\gamma(n)$ into smaller problems of determining elements of $\mathcal{L}_{\gamma,\gamma k}(x_k^*)$ for all $k$ as stated in the next theorem.

**Theorem 3.1.** *For an independent innovations tree, let there be given one leaf set $L^{(k)}$ belonging to $\mathcal{L}_{\gamma,\gamma k}(x_k^*)$ for all $k$. Then $\bigcup_{k=1}^{P_\gamma} L^{(k)} \in \mathcal{L}_\gamma(n)$. Moreover, $\mathcal{L}_{\gamma k}(n) = \mathcal{L}_{\gamma k,\gamma k}(n) = \mathcal{L}_{\gamma,\gamma k}(n)$. Also $\mu_{\gamma,\gamma k}(n)$ is a positive, non-decreasing, and discrete-concave function of $n$, $\forall k, \gamma$.*

Theorem 3.1 gives us a two step procedure to obtain the best set of $n$ leaves in the tree of $\gamma$ to estimate $V_\gamma$. We first obtain the best set of $x_k^*$ leaves in the tree of $\gamma k$ to estimate $V_{\gamma k}$ for all children $\gamma k$ of $\gamma$. We then take the union of these sets of leaves to obtain the required optimal set.

By sub-dividing the problem of obtaining optimal leaf nodes into smaller sub-problems we arrive at the following recursive technique to construct $L \in \mathcal{L}_\gamma(n)$. Starting at $\gamma$ we move downward determining how many of the $n$ leaf nodes of $L \in \mathcal{L}_\gamma(n)$ lie in the trees of the different descendants of $\gamma$ until we reach the bottom. Assume for the moment that the functions $\mu_{\gamma,\gamma k}(n)$, for all $\gamma$, are given.

**Scale-Recursive Water-filling scheme $\gamma \to \gamma k$**

**Step (a):** Split $n$ leaf nodes between the trees of $\gamma k$, $k = 1, 2, \ldots, P_\gamma$.
First determine how to split the $n$ leaf nodes between the trees of $\gamma k$ by maximizing $\sum_{k=1}^{P_\gamma} \mu_{\gamma,\gamma k}(x_k)$ over $X \in \Delta_n(\mathcal{N}_{\gamma 1}, \ldots, \mathcal{N}_{\gamma P_\gamma})$ (see (3.15)). The split is given by $X^*$ which is easily obtained using the water-filling procedure for discrete-concave functions (defined in (3.4)) since $\mu_{\gamma,\gamma k}(n)$ is discrete-concave for all $k$. Determine $L^{(k)} \in \mathcal{L}_{\gamma,\gamma k}(x_k^*)$ since $L = \bigcup_{k=1}^{P_\gamma} L^{(k)} \in \mathcal{L}_\gamma(n)$.

**Step (b):** Split $x_k^*$ nodes between the trees of child nodes of $\gamma k$.
It turns out that $L^{(k)} \in \mathcal{L}_{\gamma,\gamma k}(x_k^*)$ if and only if $L^{(k)} \in \mathcal{L}_{\gamma k}(x_k^*)$. Thus repeat Step (a) with $\gamma = \gamma k$ and $n = x_k^*$ to construct $L^{(k)}$. Stop when we have reached the bottom of the tree.

We outline an efficient implementation of the scale-recursive water-filling algorithm. This implementation first computes $L \in \mathcal{L}_\gamma(n)$ for $n = 1$ and then inductively obtains the same for larger values of $n$. Given $L \in \mathcal{L}_\gamma(n)$ we obtain $L \in \mathcal{L}_\gamma(n+1)$ as follows. Note from Step (a) above that we determine how to split the $n$ leaves at $\gamma$. We are now required to split $n+1$ leaves at $\gamma$. We easily obtain this from the earlier split of $n$ leaves using (3.4). The water-filling technique maintains the split of $n$ leaf nodes at $\gamma$ while adding just one leaf node to the tree of one of the child nodes (say $\gamma k'$) of $\gamma$. We thus have to perform Step (b) only for $k = k'$. In this way the new leaf node "percolates" down the tree until we find its location at the bottom of the tree. The pseudo-code for determining $L \in \mathcal{L}_\gamma(n)$ given $\text{var}(W_\gamma)$ for all $\gamma$ as well as the proof that the recursive water-filling algorithm can be computed in polynomial-time are available online (Ribeiro et al. [14]).

### *3.3. Uniform leaf nodes are optimal for scale-invariant trees*

The symmetry in scale-invariant trees forces the optimal solution to take a particular form irrespective of the variances of the innovations $W_\gamma$. We use the following notion of *uniform split* to prove that in a scale-invariant tree a more or less equal spread of sample leaf nodes across the tree gives the best linear estimate of the root.

**Definition 3.2.** Given a scale-invariant tree, a vector of leaf nodes $L$ has uniform split of size $n$ at node $\gamma$ if $|L_\gamma| = n$ and $|L_{\gamma k}|$ is either $\lfloor \frac{n}{P_\gamma} \rfloor$ or $\lfloor \frac{n}{P_\gamma} \rfloor + 1$ for all values of $k$. It follows that $\#\{k : |L_{\gamma k}| = \lfloor \frac{n}{P_\gamma} \rfloor + 1\} = n - P_\gamma \lfloor \frac{n}{P_\gamma} \rfloor$.

**Definition 3.3.** Given a scale-invariant tree, a vector of leaf nodes is called a uniform leaf sample if it has a uniform split at all tree nodes.

The next theorem gives the optimal leaf node set for scale-invariant trees.

**Theorem 3.2.** *Given a scale-invariant tree, the uniform leaf sample of size $n$ gives the best LMMSE estimate of the tree-root among all possible choices of $n$ leaf nodes.*

*Proof.* For a scale-invariant tree, $\mu_{\gamma,\gamma k}(n)$ is identical for all $k$ given any location $\gamma$. Corollary 3.1 and Theorem 3.1 then prove the theorem. □

## 4. Covariance trees

In this section we prove optimal and worst case solutions for covariance trees. For the optimal solutions we leverage our results for independent innovations trees and for the worst case solutions we employ eigenanalysis. We begin by formulating the problem.

### *4.1. Problem formulation*

Let us compute the LMMSE of estimating the root $V_\emptyset$ given a set of leaf nodes $L$ of size $n$. Recall that for a covariance tree the correlation between any leaf node and the root node is identical. We denote this correlation by $\rho$. Denote an $i \times j$ matrix with all elements equal to 1 by $\mathbf{1}_{i \times j}$. It is well known (Stark and Woods [17]) that the optimal linear estimate of $V_\emptyset$ given $L$ (assuming zero-mean random variables) is given by $\rho \mathbf{1}_{1 \times n} Q_L^{-1} L$, where $Q_L$ is the covariance matrix of $L$ and that the resulting LMMSE is

$$\begin{aligned}(4.1) \quad \mathcal{E}(V_\emptyset | L) &= \mathrm{var}(V_\emptyset) - \mathrm{cov}(L, V_\emptyset)^T Q_L^{-1} \mathrm{cov}(L, V_\emptyset) \\ &= \mathrm{var}(V_\emptyset) - \rho^2 \mathbf{1}_{1 \times n} Q_L^{-1} \mathbf{1}_{n \times 1}.\end{aligned}$$

Clearly obtaining the best and worst-case choices for $L$ is equivalent to maximizing and minimizing the sum of the elements of $Q_L^{-1}$. The exact value of $\rho$ does not affect the solution. We assume that no element of $L$ can be expressed as a linear combination of the other elements of $L$ which implies that $Q_L$ is invertible.

### *4.2. Optimal solutions*

We use our results of Section 3 for independent innovations trees to determine the optimal solutions for covariance trees. Note from (4.2) that the estimation error for a covariance tree is a function only of the covariance *between leaf nodes*. Exploiting this fact, we first construct an independent innovations tree whose leaf nodes have the same correlation structure as that of the covariance tree and then prove that both trees must have the same optimal solution. Previous results then provide the optimal solution for the independent innovations tree which is also optimal for the covariance tree.

**Definition 4.1.** A matched innovations tree of a given covariance tree with positive correlation progression is an independent innovations tree with the following properties. It has (1) the same topology (2) and the same correlation structure between leaf nodes as the covariance tree, and (3) the root is equally correlated with all leaf nodes (though the exact value of the correlation between the root and a leaf node may differ from that of the covariance tree).

All covariance trees with positive correlation progression have corresponding matched innovations trees. We construct a matched innovations tree for a given covariance tree as follows. Consider an independent innovations tree with the same topology as the covariance tree. Set $\varrho_\gamma = 1$ for all $\gamma$,

$$\text{(4.2)} \qquad \text{var}(V_\emptyset) = c_0$$

and

$$\text{(4.3)} \qquad \text{var}(W^{(j)}) = c_j - c_{j-1}, \ j = 1, 2, \ldots, D,$$

where $c_j$ is the covariance of leaf nodes of the covariance tree with proximity $j$ and $\text{var}(W^{(j)})$ is the common variance of all innovations of the independent innovations tree at scale $j$. Call $c'_j$ the covariance of leaf nodes with proximity $j$ in the independent innovations tree. From (2.1) we have

$$\text{(4.4)} \qquad c'_j = \text{var}(V_\emptyset) + \sum_{k=1}^{j} \text{var}\left(W^{(k)}\right), \ j = 1, \ldots, D.$$

Thus, $c'_j = c_j$ for all $j$ and hence this independent innovations tree is the required matched innovations tree.

The next lemma relates the optimal solutions of a covariance tree and its matched innovations tree.

**Lemma 4.1.** *A covariance tree with positive correlation progression and its matched innovations tree have the same optimal leaf sets.*

*Proof.* Note that (4.2) applies to any tree whose root is equally correlated with all its leaves. This includes both the covariance tree and its matched innovations tree. From (4.2) we see that the choice of $L$ that maximizes the sum of elements of $Q_L^{-1}$ is optimal. Since $Q_L^{-1}$ is identical for both the covariance tree and its matched innovations tree for any choice of $L$, they must have the same optimal solution. □

For a symmetric covariance tree that has positive correlation progression, the optimal solution takes on a specific form irrespective of the actual covariance between leaf nodes.

**Theorem 4.1.** *Given a symmetric covariance tree that has positive correlation progression, the uniform leaf sample of size $n$ gives the best LMMSE of the tree-root among all possible choices of $n$ leaf nodes.*

*Proof.* Form a matched innovations tree using the procedure outlined previously. This tree is by construction a scale-invariant tree. The result then follows from Theorem 3.2 and Lemma 4.1. □

While the uniform leaf sample is the optimal solution for a symmetric covariance tree with positive correlation progression, it is surprisingly the worst case solution for certain trees with a different correlation structure, which we prove next.

### *4.3. Worst case solutions*

The *worst case solution* is any choice of $L \in \Lambda_\emptyset(n)$ that maximizes $\mathcal{E}(V_\emptyset|L)$. We now highlight the fact that the best and worst case solutions can change dramatically depending on the correlation structure of the tree. Of particular relevance to our discussion is the set of *clustered leaf nodes* defined as follows.

**Definition 4.2.** The set consisting of all leaf nodes of the tree of $V_\gamma$ is called the set of clustered leaves of $\gamma$.

We provide the worst case solutions for covariance trees in which every node (with the exception of the leaves) has the same number of child nodes. The following theorem summarizes our results.

**Theorem 4.2.** *Consider a covariance tree of depth $D$ in which every node (excluding the leaves) has the same number of child nodes $\sigma$. Then for leaf sets of size $\sigma^p$, $p = 0, 1, \ldots, D$, the worst case solution when the tree has positive correlation progression is given by the sets of clustered leaves of $\gamma$, where $\gamma$ is any node at scale $D - p$. The worst case solution is given by the sets of uniform leaf nodes when the tree has negative correlation progression.*

Theorem 4.2 gives us the intuition that "more correlated" leaf nodes give worse estimates of the root. In the case of covariance trees with positive correlation progression, clustered leaf nodes are strongly correlated when compared to uniform leaf nodes. The opposite is true in the negative correlation progression case. Essentially if leaf nodes are highly correlated then they contain more redundant information which leads to poor estimation of the root.

While we have proved the optimal solution for covariance trees with positive correlation progression. we have not yet proved the same for those with negative correlation progression. Based on the intuition just gained we make the following conjecture.

**Conjecture 4.1.** *Consider a covariance tree of depth $D$ in which every node (excluding the leaves) has the same number of child nodes $\sigma$. Then for leaf sets of size $\sigma^p$, $p = 0, 1, \ldots, D$, the optimal solution when the tree has negative correlation progression is given by the sets of clustered leaves of $\gamma$, where $\gamma$ is any node at scale $D - p$.*

Using numerical techniques we support this conjecture in the next section.

## 5. Numerical results

In this section, using the scale-recursive water-filling algorithm we evaluate the optimal leaf sets for independent innovations trees that are not scale-invariant. In addition we provide numerical support for Conjecture 4.1.

### *5.1. Independent innovations trees: scale-recursive water-filling*

We consider trees with depth $D = 3$ and in which all nodes have at most two child nodes. The results demonstrate that the optimal leaf sets are a function of the correlation structure and topology of the multiscale trees.

In Fig. 4(a) we plot the optimal leaf node sets of different sizes for a scale-invariant tree. As expected the uniform leaf nodes sets are optimal.

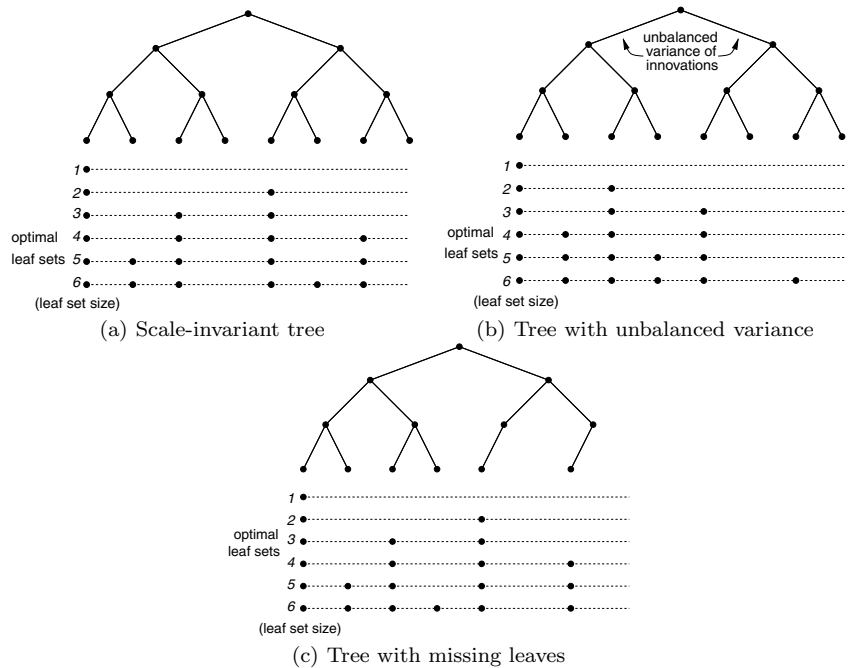

(a) Scale-invariant tree  (b) Tree with unbalanced variance

(c) Tree with missing leaves

FIG 4. *Optimal leaf node sets for three different independent innovations trees:* (a) *scale-invariant tree,* (b) *symmetric tree with unbalanced variance of innovations at scale* 1, *and* (c) *tree with missing leaves at the finest scale. Observe that the uniform leaf node sets are optimal in* (a) *as expected. In* (b), *however, the nodes on the left half of the tree are more preferable to those on the right. In* (c) *the solution is similar to* (a) *for optimal sets of size* $n = 5$ *or lower but changes for* $n = 6$ *due to the missing nodes.*

We consider a symmetric tree in Fig. 4(b), that is a tree in which all nodes have the same number of children (excepting leaf nodes). All parameters are constant within each scale except for the variance of the innovations $W_\gamma$ at scale 1. The variance of the innovation on the right side is five times larger than the variance of the innovation on the left. Observe that leaves on the left of the tree are now preferable to those on the right and hence dominate the optimal sets. Comparing this result to Fig. 4(a) we see that the optimal sets are dependent on the correlation structure of the tree.

In Fig. 4(c) we consider the same tree as in Fig. 4(a) with two leaf nodes missing. These two leaves do not belong to the optimal leaf sets of size $n = 1$ to $n = 5$ in Fig. 4(a) but are elements of the optimal set for $n = 6$. As a result the optimal sets of size 1 to 5 in Fig. 4(c) are identical to those in Fig. 4(a) whereas that for $n = 6$ differs. This result suggests that the optimal sets depend on the tree topology.

Our results have important implications for applications because situations arise where we must model physical processes using trees with different correlation structures and topologies. For example, if the process to be measured is non-stationary over space then the multiscale tree may be unbalanced as in Fig. 4(b). In some applications it may not be possible to sample at certain locations due to physical constraints. We would thus have to exclude certain leaf nodes in our analysis as in Fig. 4(c).

The above experiments with tree-depth $D = 3$ are "toy-examples" to illustrate key concepts. In practice, the water-filling algorithm can solve much larger real-world problems with ease. For example, on a Pentium IV machine running Matlab, the water-filling algorithm takes 22 seconds to obtain the optimal leaf set of size 100 to estimate the root of a binary tree with depth 11, that is a tree with 2048 leaves.

## *5.2. Covariance trees: best and worst cases*

This section provides numerical support for Conjecture 4.1 that states that the clustered leaf node sets are optimal for covariance trees with negative correlation progression. We employ the WIG tree, a covariance tree in which each node has $\sigma = 2$ child nodes (Ma and Ji [12]). We provide numerical support for our claim using a WIG model of depth $D = 6$ possessing a fractional Gaussian noise-like[2] correlation structure corresponding to $H = 0.8$ and $H = 0.3$. To be precise, we choose the WIG model parameters such that the variance of nodes at scale $j$ is proportional to $2^{-2jH}$ (see Ma and Ji [12] for further details). Note that $H > 0.5$ corresponds to positive correlation progression while $H \leq 0.5$ corresponds to negative correlation progression.

Fig. 5 compares the LMMSE of the estimated root node (normalized by the variance of the root) of the uniform and clustered sampling patterns. Since an exhaustive search of all possible patterns is very computationally expensive (for example there are over $10^{18}$ ways of choosing 32 leaf nodes from among 64) we instead compute the LMMSE for $10^4$ randomly selected patterns. Observe that the clustered pattern gives the smallest LMMSE for the tree with negative correlation progression in Fig. 5(a) supporting our Conjecture 4.1 while the uniform pattern gives the smallest LMMSE for the positively correlation progression one in Fig. 5(b) as stated in Theorem 4.1. As proved in Theorem 4.2, the clustered and uniform patterns give the worst LMMSE for the positive and negative correlation progression cases respectively.

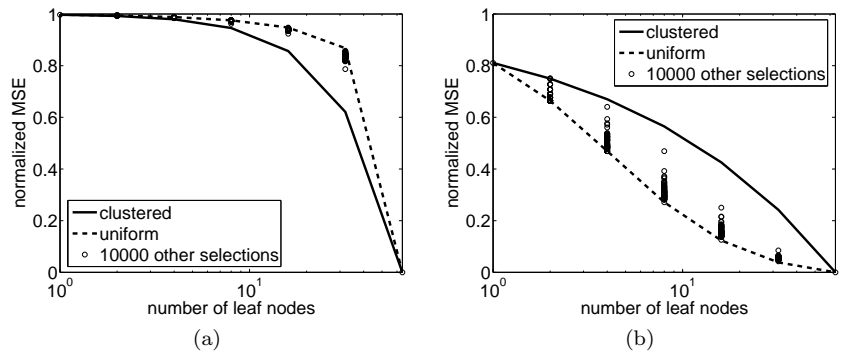

FIG 5. *Comparison of sampling schemes for a WIG model with* (a) *negative correlation progression and* (b) *positive correlation progression. Observe that the clustered nodes are optimal in* (a) *while the uniform is optimal in* (b). *The uniform and the clustered leaf sets give the worst performance in* (a) *and* (b) *respectively, as expected from our theoretical results.*

---

[2]Fractional Gaussian noise is the increments process of fBm (Mandelbrot and Ness [13]).

## 6. Related work

Earlier work has studied the problem of designing optimal samples of size $n$ to linearly estimate the sum total of a process. For a one dimensional process which is wide-sense stationary with positive and convex correlation, within a class of unbiased estimators of the sum of the population, it was shown that systematic sampling of the process (uniform patterns with random starting points) is optimal (Hájek [6]).

For a two dimensional process on an $n_1 \times n_2$ grid with positive and convex correlation it was shown that an optimal sampling scheme does not lie in the class of schemes that ensure equal inclusion probability of $n/(n_1 n_2)$ for every point on the grid (Bellhouse [2]). In Bellhouse [2], an "optimal scheme" refers to a sampling scheme that achieves a particular lower bound on the error variance. The requirement of equal inclusion probability guarantees an unbiased estimator. The optimal schemes *within* certain sub-classes of this larger "equal inclusion probability" class were obtained using systematic sampling. More recent analysis refines these results to show that optimal designs do exist in the equal inclusion probability class for certain values of $n$, $n_1$, and $n_2$ and are obtained by Latin square sampling (Lawry and Bellhouse [10], Salehi [16]).

Our results differ from the above works in that we provide optimal solutions for the entire class of linear estimators and study a different set of random processes.

Other work on sampling fractional Brownian motion to estimate its Hurst parameter demonstrated that geometric sampling is superior to uniform sampling (Vidàcs and Virtamo [18]).

Recent work compared different probing schemes for traffic estimation through numerical simulations (He and Hou [7]). It was shown that a scheme which used uniformly spaced probes outperformed other schemes that used clustered probes. These results are similar to our findings for independent innovation trees and covariance trees with positive correlation progression.

## 7. Conclusions

This paper has addressed the problem of obtaining optimal leaf sets to estimate the root node of two types of multiscale stochastic processes: independent innovations trees and covariance trees. Our findings are particularly useful for applications which require the estimation of the sum total of a correlated population from a finite sample.

We have proved for an independent innovations tree that the optimal solution can be obtained using an efficient water-filling algorithm. Our results show that the optimal solutions can vary drastically depending on the correlation structure of the tree. For covariance trees with positive correlation progression as well as scale-invariant trees we obtained that uniformly spaced leaf nodes are optimal. However, uniform leaf nodes give the worst estimates for covariance trees with negative correlation progression. Numerical experiments support our conjecture that clustered nodes provide the optimal solution for covariance trees with negative correlation progression.

This paper raises several interesting questions for future research. The general problem of determining which $n$ random variables from a given set provide the best linear estimate of another random variable that is not in the same set is an NP-hard problem. We, however, devised a fast polynomial-time algorithm to solve one problem of this type, namely determining the optimal leaf set for an independent innovations tree. Clearly, the structure of independent innovations trees was an important factor that enabled a fast algorithm. The question arises as to whether there are similar problems that have polynomial-time solutions.

We have proved optimal results for covariance trees by reducing the problem to one for independent innovations trees. Such techniques of reducing one optimization problem to another problem that has an efficient solution can be very powerful. If a problem can be reduced to one of determining optimal leaf sets for independent innovations trees in polynomial-time, then its solution is also polynomial-time. Which other problems are malleable to this reduction is an open question.

**Appendix**

*Proof of Lemma 3.1.* We first prove the following statement.
*Claim* (1): If there exists $X^* = [x_k^*] \in \Delta_n(M_1, \ldots, M_P)$ that has the following property:

$$\psi_i(x_i^*) - \psi_i(x_i^* - 1) \geq \psi_j(x_j^* + 1) - \psi_j(x_j^*), \tag{7.1}$$

$\forall i \neq j$ such that $x_i^* > 0$ and $x_j^* < M_j$, then

$$h(n) = \sum_{k=1}^{P} \psi_k(x_k^*). \tag{7.2}$$

We then prove that such an $X^*$ always exists and can be constructed using the water-filling technique.

Consider any $\widehat{X} \in \Delta_n(M_1, \ldots, M_P)$. Using the following steps, we transform the vector $\widehat{X}$ two elements at a time to obtain $X^*$.

Step 1: (Initialization) Set $X = \widehat{X}$.
Step 2: If $X \neq X^*$, then since the elements of both $X$ and $X^*$ sum up to $n$, there must exist a pair $i, j$ such that $x_i \neq x_i^*$ and $x_j \neq x_j^*$. Without loss of generality assume that $x_i < x_i^*$ and $x_j > x_j^*$. This assumption implies that $x_i^* > 0$ and $x_j^* < M_j$. Now form vector $Y$ such that $y_i = x_i + 1$, $y_j = x_j - 1$, and $y_k = x_k$ for $k \neq i, j$. From (7.1) and the concavity of $\psi_i$ and $\psi_j$ we have

$$\begin{aligned}
\psi_i(y_i) - \psi_i(x_i) &= \psi_i(x_i + 1) - \psi_i(x_i) \geq \psi_i(x_i^*) - \psi_i(x_i^* - 1) \\
&\geq \psi_j(x_j^* + 1) - \psi_j(x_j^*) \geq \psi_j(x_j) - \psi_j(x_j - 1) \\
&\geq \psi_j(x_j) - \psi_j(y_j).
\end{aligned} \tag{7.3}$$

As a consequence

$$\sum_k (\psi_k(y_k) - \psi_k(x_k)) = \psi_i(y_i) - \psi_i(x_i) + \psi_j(y_j) - \psi_j(x_j) \geq 0. \tag{7.4}$$

Step 3: If $Y \neq X^*$ then set $X = Y$ and repeat Step 2, otherwise stop.

After performing the above steps at most $\sum_k M_k$ times, $Y = X^*$ and (7.4) gives

$$\sum_k \psi_k(x_k^*) = \sum_k \psi_k(y_k) \geq \sum_k \psi_k(\widehat{x}_k). \tag{7.5}$$

This proves Claim (1).

Indeed for any $\widetilde{X} \neq X^*$ satisfying (7.1) we must have $\sum_k \psi_k(\widetilde{x}_k) = \sum_k \psi_k(x_k^*)$. We now prove the following claim by induction.

*Claim* (2): $G^{(n)} \in \Delta_n(M_1, \ldots, M_P)$ and $G^{(n)}$ satisfies (7.1).
(Initial Condition) The claim is trivial for $n = 0$.
(Induction Step) Clearly from (3.4) and (3.5)

$$\text{(7.6)} \qquad \sum_k g_k^{(n+1)} = 1 + \sum_k g_k^{(n)} = n + 1,$$

and $0 \leq g_k^{(n+1)} \leq M_k$. Thus $G^{(n+1)} \in \Delta_{n+1}(M_1, \ldots, M_P)$. We now prove that $G^{(n+1)}$ satisfies property (7.1). We need to consider pairs $i, j$ as in (7.1) for which either $i = m$ or $j = m$ because all other cases directly follow from the fact that $G^{(n)}$ satisfies (7.1).

Case (i) $j = m$, where $m$ is defined as in (3.5). Assuming that $g_m^{(n+1)} < M_m$, for all $i \neq m$ such that $g_i^{(n+1)} > 0$ we have

$$\text{(7.7)} \qquad \begin{aligned} \psi_i\left(g_i^{(n+1)}\right) - \psi_i\left(g_i^{(n+1)} - 1\right) &= \psi_i\left(g_i^{(n)}\right) - \psi_i\left(g_i^{(n)} - 1\right) \\ &\geq \psi_m\left(g_m^{(n)} + 1\right) - \psi_m\left(g_m^{(n)}\right) \\ &\geq \psi_m\left(g_m^{(n)} + 2\right) - \psi_m\left(g_m^{(n)} + 1\right) \\ &= \psi_m\left(g_m^{(n+1)} + 1\right) - \psi_m\left(g_m^{(n+1)}\right). \end{aligned}$$

Case (ii) $i = m$. Consider $j \neq m$ such that $g_j^{(n+1)} < M_j$. We have from (3.5) that

$$\text{(7.8)} \qquad \begin{aligned} \psi_m\left(g_m^{(n+1)}\right) - \psi_m\left(g_m^{(n+1)} - 1\right) &= \psi_m\left(g_m^{(n)} + 1\right) - \psi_m\left(g_m^{(n)}\right) \\ &\geq \psi_j\left(g_j^{(n)} + 1\right) - \psi_j\left(g_j^{(n)}\right) \\ &= \psi_j\left(g_j^{(n+1)} + 1\right) - \psi_j\left(g_j^{(n+1)}\right). \end{aligned}$$

Thus Claim (2) is proved.

It only remains to prove the next claim.

*Claim* (3): $h(n)$, or equivalently $\sum_k \psi_k(g_k^{(n)})$, is non-decreasing and discrete-concave.

Since $\psi_k$ is non-decreasing for all $k$, from (3.4) we have that $\sum_k \psi_k(g_k^{(n)})$ is a non-decreasing function of $n$. We have from (3.5)

$$\text{(7.9)} \qquad \begin{aligned} h(n+1) - h(n) &= \sum_k \left(\psi_k(g_k^{(n+1)}) - \psi_k(g_k^{(n)})\right) \\ &= \max_{k: g_k^{(n)} < M_k} \left\{\psi_k(g_k^{(n)} + 1) - \psi_k(g_k^{(n)})\right\}. \end{aligned}$$

From the concavity of $\psi_k$ and the fact that $g_k^{(n+1)} \geq g_k^{(n)}$ we have that

$$\text{(7.10)} \qquad \psi_k(g_k^{(n)} + 1) - \psi_k(g_k^{(n)}) \geq \psi_k(g_k^{(n+1)} + 1) - \psi_k(g_k^{(n+1)}),$$

for all $k$. Thus from (7.10) and (7.10), $h(n)$ is discrete-concave. □

*Proof of Corollary 3.1.* Set $x_k^* = \lfloor \frac{n}{P} \rfloor$ for $1 \leq k \leq P - n + P\lfloor \frac{n}{P} \rfloor$ and $x_k^* = 1 + \lfloor \frac{n}{P} \rfloor$ for all other $k$. Then $X^* = [x_k^*] \in \Delta_n(M_1, \ldots, M_P)$ and $X^*$ satisfies (7.1) from which the result follows. □

The following two lemmas are required to prove Theorem 3.1.

**Lemma 7.1.** *Given independent random variables $A, W, F$, define $Z$ and $E$ through $Z := \zeta A + W$ and $E := \eta Z + F$ where $\zeta, \eta$ are constants. We then have the result*

$$\text{(7.11)} \qquad \frac{\text{var}(A)}{\text{cov}(A,E)^2} \cdot \frac{\text{cov}(Z,E)^2}{\text{var}(Z)} = \frac{\zeta^2 + \text{var}(W)/\text{var}(A)}{\zeta^2} \geq 1.$$

*Proof.* Without loss of generality assume all random variables have zero mean. We have

$$\text{(7.12)} \qquad \text{cov}(E, Z) = \mathbb{E}(\eta Z^2 + FZ) = \eta \text{var}(Z),$$

$$\text{(7.13)} \qquad \text{cov}(A, E) = \mathbb{E}((\eta(\zeta A + W) + F)A)\zeta\eta\text{var}(A),$$

and

$$\text{(7.14)} \qquad \text{var}(Z) = \mathbb{E}(\zeta^2 A^2 + W^2 + 2\zeta AW) = \zeta^2 \text{var}(A) + \text{var}(W).$$

Thus from (7.12), (7.13) and (7.14)

$$\text{(7.15)} \qquad \frac{\text{cov}(Z,E)^2}{\text{var}(Z)} \cdot \frac{\text{var}(A)}{\text{cov}(A,E)^2} = \frac{\eta^2 \text{var}(Z)}{\zeta^2 \eta^2 \text{var}(A)} = \frac{\zeta^2 + \text{var}(W)/\text{var}(A)}{\zeta^2} \geq 1.$$

□

**Lemma 7.2.** *Given a positive function $z_i, i \in \mathbb{Z}$ and constant $\alpha > 0$ such that*

$$\text{(7.16)} \qquad r_i := \frac{1}{1 - \alpha z_i}$$

*is positive, discrete-concave, and non-decreasing, we have that*

$$\text{(7.17)} \qquad \delta_i := \frac{1}{1 - \beta z_i}$$

*is also positive, discrete-concave, and non-decreasing for all $\beta$ with $0 < \beta \leq \alpha$.*

*Proof.* Define $\kappa_i := z_i - z_{i-1}$. Since $z_i$ is positive and $r_i$ is positive and non-decreasing, $\alpha z_i < 1$ and $z_i$ must increase with $i$, that is $\kappa_i \geq 0$. This combined with the fact that $\beta z_i \leq \alpha z_i < 1$ guarantees that $\delta_i$ must be positive and non-decreasing.

It only remains to prove the concavity of $\delta_i$. From (7.16)

$$\text{(7.18)} \qquad r_{i+1} - r_i = \frac{\alpha(z_{i+1} - z_i)}{(1 - \alpha z_{i+1})(1 - \alpha z_i)} = \alpha \kappa_{i+1} r_{i+1} r_i.$$

We are given that $r_i$ is discrete-concave, that is

$$\text{(7.19)} \qquad \begin{aligned} 0 &\geq (r_{i+2} - r_{i+1}) - (r_{i+1} - r_i) \\ &= \alpha r_i r_{i+1} \left[ \kappa_{i+2} \left( \frac{1 - \alpha z_i}{1 - \alpha z_{i+2}} \right) - \kappa_{i+1} \right]. \end{aligned}$$

Since $r_i > 0 \ \forall i$, we must have

$$\text{(7.20)} \qquad \left[ \kappa_{i+2} \left( \frac{1 - \alpha z_i}{1 - \alpha z_{i+2}} \right) - \kappa_{i+1} \right] \leq 0.$$

Similar to (7.20) we have that

$$\text{(7.21)} \qquad (\delta_{i+2} - \delta_{i+1}) - (\delta_{i+1} - \delta_i) = \beta \delta_i \delta_{i+1} \left[ \kappa_{i+2} \left( \frac{1 - \beta z_i}{1 - \beta z_{i+2}} \right) - \kappa_{i+1} \right].$$

Since $\delta_i > 0 \ \forall i$, for the concavity of $\delta_i$ it suffices to show that

$$\left[\kappa_{i+2} \frac{1-\beta z_i}{1-\beta z_{i+2}} - \kappa_{i+1}\right] \leq 0. \tag{7.22}$$

Now

$$\frac{1-\alpha z_i}{1-\alpha z_{i+2}} - \frac{1-\beta z_i}{1-\beta z_{i+2}} = \frac{(\alpha-\beta)(z_{i+2}-z_i)}{(1-\alpha z_{i+2})(1-\beta z_{i+2})} \geq 0. \tag{7.23}$$

Then (7.20) and (7.23) combined with the fact that $\kappa_i \geq 0$, $\forall i$ proves (7.22). □

*Proof of Theorem 3.1.* We split the theorem into three claims.

*Claim (1):* $L^* := \cup_k L^{(k)}(x_k^*) \in \mathcal{L}_\gamma(n)$.

From (3.10), (3.11), and (3.13) we obtain

$$\begin{aligned}
\mu_\gamma(n) + \frac{P_\gamma - 1}{\mathrm{var}(V_\gamma)} &= \max_{L \in \Lambda_\gamma(n)} \sum_{k=1}^{P_\gamma} \mathcal{E}(V_\gamma | L_{\gamma k})^{-1} \\
&\leq \max_{X \in \Delta_n(\mathcal{N}_{\gamma 1},\ldots,\mathcal{N}_{\gamma P_\gamma})} \sum_{k=1}^{P_\gamma} \mu_{\gamma,\gamma k}(x_k).
\end{aligned} \tag{7.24}$$

Clearly $L^* \in \Lambda_\gamma(n)$. We then have from (3.10) and (3.11) that

$$\begin{aligned}
\mu_\gamma(n) + \frac{P_\gamma - 1}{\mathrm{var}(V_\gamma)} &\geq \mathcal{E}(V_\gamma|L^*)^{-1} + \frac{P_\gamma - 1}{\mathrm{var}(V_\gamma)} = \sum_{k=1}^{P_\gamma} \mathcal{E}(V_\gamma|L_{\gamma k}^*)^{-1} \\
&= \sum_{k=1}^{P_\gamma} \mu_{\gamma,\gamma k}(x_k^*) = \max_{X \in \Delta_n(\mathcal{N}_{\gamma 1},\ldots,\mathcal{N}_{\gamma P_\gamma})} \sum_{k=1}^{P_\gamma} \mu_{\gamma,\gamma k}(x_k).
\end{aligned} \tag{7.25}$$

Thus from (7.25) and (7.26) we have

$$\mu_\gamma(n) = \mathcal{E}(V_\gamma|L^*)^{-1} = \max_{X \in \Delta_n(\mathcal{N}_{\gamma 1},\ldots,\mathcal{N}_{\gamma P_\gamma})} \sum_{k=1}^{P_\gamma} \mu_{\gamma,\gamma k}(x_k) - \frac{P_\gamma - 1}{\mathrm{var}(V_\gamma)}, \tag{7.26}$$

which proves Claim (1).

*Claim (2):* If $L \in \mathcal{L}_{\gamma k}(n)$ then $L \in \mathcal{L}_{\gamma,\gamma k}(n)$ and vice versa.
Denote an arbitrary leaf node of the tree of $\gamma k$ as $E$. Then $V_\gamma$, $V_{\gamma k}$, and $E$ are related through

$$V_{\gamma k} = \varrho_{\gamma k} V_\gamma + W_{\gamma k}, \tag{7.27}$$

and

$$E = \eta V_{\gamma k} + F \tag{7.28}$$

where $\eta$ and $\varrho_{\gamma k}$ are scalars and $W_{\gamma k}$, $F$ and $V_\gamma$ are independent random variables. We note that by definition $\mathrm{var}(V_\gamma) > 0 \ \forall \gamma$ (see Definition 2.5). From Lemma 7.1

we have

$$(7.29) \quad \frac{\text{cov}(V_{\gamma k}, E)}{\text{cov}(V_\gamma, E)} = \left(\frac{\text{var}(V_{\gamma k})}{\text{var}(V_\gamma)}\right)^{1/2} \left(\frac{\varrho_{\gamma k}^2 + \frac{\text{var}(W_{\gamma k})}{\text{var}(V_\gamma)}}{\varrho_{\gamma k}^2}\right)^{1/2}$$

$$=: \xi_{\gamma,k} \geq \left(\frac{\text{var}(V_{\gamma k})}{\text{var}(V_\gamma)}\right)^{1/2}.$$

From (7.30) we see that $\xi_{\gamma,k}$ is not a function of $E$.

Denote the covariance between $V_\gamma$ and leaf node vector $L = [\ell_i] \in \Lambda_{\gamma k}(n)$ as $\Theta_{\gamma,L} = [\text{cov}(V_\gamma, \ell_i)]^T$. Then (7.30) gives

$$(7.30) \quad \Theta_{\gamma k, L} = \xi_{\gamma,k} \Theta_{\gamma,L}.$$

From (4.2) we have

$$(7.31) \quad \mathcal{E}(V_\gamma | L) = \text{var}(V_\gamma) - \varphi(\gamma, L)$$

where $\varphi(\gamma, L) = \Theta_{\gamma,L}^T Q_L^{-1} \Theta_{\gamma,L}$. Note that $\varphi(\gamma, L) \geq 0$ since $Q_L^{-1}$ is positive semi-definite. Using (7.30) we similarly get

$$(7.32) \quad \mathcal{E}(V_{\gamma k} | L) = \text{var}(V_{\gamma k}) - \frac{\varphi(\gamma, L)}{\xi_{\gamma,k}^2}.$$

From (7.31) and (7.32) we see that $\mathcal{E}(V_\gamma | L)$ and $\mathcal{E}(V_{\gamma k} | L)$ are both minimized over $L \in \Lambda_{\gamma k}(n)$ by the same leaf vector that maximizes $\varphi(\gamma, L)$. This proves Claim (2).

*Claim* (3): $\mu_{\gamma,\gamma k}(n)$ is a positive, non-decreasing, and discrete-concave function of $n$, $\forall k, \gamma$.

We start at a node $\gamma$ at one scale from the bottom of the tree and then move up the tree.

*Initial Condition:* Note that $V_{\gamma k}$ is a leaf node. From (2.1) and (**??**) we obtain

$$(7.33) \quad \mathcal{E}(V_\gamma | V_{\gamma k}) = \text{var}(V_\gamma) - \frac{(\varrho_{\gamma k} \text{var}(V_\gamma))^2}{\text{var}(V_{\gamma k})} \leq \text{var}(V_\gamma).$$

For our choice of $\gamma$, $\mu_{\gamma,\gamma k}(1)$ corresponds to $\mathcal{E}(V_\gamma | V_{\gamma k})^{-1}$ and $\mu_{\gamma,\gamma k}(0)$ corresponds to $1/\text{var}(V_\gamma)$. Thus from (7.33), $\mu_{\gamma,\gamma k}(n)$ is positive, non-decreasing, and discrete-concave (trivially since $n$ takes only two values here).

*Induction Step:* Given that $\mu_{\gamma,\gamma k}(n)$ is a positive, non-decreasing, and discrete-concave function of $n$ for $k = 1, \ldots, P_\gamma$, we prove the same when $\gamma$ is replaced by $\gamma \uparrow$. Without loss of generality choose $k$ such that $(\gamma \uparrow)k = \gamma$. From (3.11), (3.13), (7.31), (7.32) and Claim (2), we have for $L \in \mathcal{L}_\gamma(n)$

$$(7.34) \quad \begin{aligned} \mu_\gamma(n) &= \frac{1}{\text{var}(V_\gamma)} \cdot \frac{1}{1 - \frac{\varphi(\gamma,L)}{\text{var}(V_\gamma)}}, \quad \text{and} \\ \mu_{\gamma\uparrow,k}(n) &= \frac{1}{\text{var}(V_{\gamma\uparrow})} \cdot \frac{1}{1 - \frac{\varphi(\gamma,L)}{\xi_{\gamma\uparrow,k}^2 \text{var}(V_{\gamma\uparrow})}}. \end{aligned}$$

From (7.26), the assumption that $\mu_{\gamma,\gamma k}(n)$ $\forall k$ is a positive, non-decreasing, and discrete-concave function of $n$, and Lemma 3.1 we have that $\mu_\gamma(n)$ is a non-decreasing and discrete-concave function of $n$. Note that by definition (see (3.11)) $\mu_\gamma(n)$ is positive. This combined with (2.1), (7.35), (7.30) and Lemma 7.2, then prove that $\mu_{\gamma\uparrow,k}(n)$ is also positive, non-decreasing, and discrete-concave. □

We now prove a lemma to be used to prove Theorem 4.2. As a first step we compute the leaf arrangements $L$ which maximize and minimize the sum of all elements of $Q_L = [q_{i,j}(L)]$. We restrict our analysis to a covariance tree with depth $D$ and in which each node (excluding leaf nodes) has $\sigma$ child nodes. We introduce some notation. Define

(7.35) $\quad \Gamma^{(u)}(p) := \{L : L \in \Lambda_\emptyset(\sigma^p) \text{ and } L \text{ is a uniform leaf node set}\}$ and

(7.36) $\quad \Gamma^{(c)}(p) := \{L : L \text{ is a clustered leaf set of a node at scale } D - p\}$

for $p = 0, 1, \ldots, D$. We number nodes at scale $m$ in an arbitrary order from $q = 0, 1, \ldots, \sigma^m - 1$ and refer to a node by the pair $(m, q)$.

**Lemma 7.3.** *Assume a positive correlation progression. Then, $\sum_{i,j} q_{i,j}(L)$ is minimized over $L \in \Lambda_\emptyset(\sigma^p)$ by every $L \in \Gamma^{(u)}(p)$ and maximized by every $L \in \Gamma^{(c)}(p)$. For a negative correlation progression, $\sum_{i,j} q_{i,j}(L)$ is maximized by every $L \in \Gamma^{(u)}(p)$ and minimized by every $L \in \Gamma^{(c)}(p)$.*

*Proof.* Set $p$ to be an arbitrary element in $\{1, \ldots, D-1\}$. The cases of $p = 0$ and $p = D$ are trivial. Let $\vartheta_m = \#\{q_{i,j}(L) \in Q_L : q_{i,j}(L) = c_m\}$ be the number of elements of $Q_L$ equal to $c_m$. Define $a_m := \sum_{k=0}^{m} \vartheta_k, m \geq 0$ and set $a_{-1} = 0$. Then

(7.37)
$$\begin{aligned}
\sum_{i,j} q_{i,j} &= \sum_{m=0}^{D} c_m \vartheta_m = \sum_{m=0}^{D-1} c_m(a_m - a_{m-1}) + c_D \vartheta_D \\
&= \sum_{m=0}^{D-1} c_m a_m - \sum_{m=-1}^{D-2} c_{m+1} a_m + c_D \vartheta_D \\
&= \sum_{m=0}^{D-2} (c_m - c_{m+1}) a_m + c_{D-1} a_{D-1} - c_0 a_{-1} + c_D \vartheta_D \\
&= \sum_{m=0}^{D-2} (c_m - c_{m+1}) a_m + \text{constant},
\end{aligned}$$

where we used the fact that $a_{D-1} = a_D - \vartheta_D$ is a constant independent of the choice of $L$, since $\vartheta_D = \sigma^p$ and $a_D = \sigma^{2p}$.

We now show that $L \in \Gamma^{(u)}(p)$ maximizes $a_m, \forall m$ while $L \in \Gamma^{(c)}(p)$ minimizes $a_m, \forall m$. First we prove the results for $L \in \Gamma^{(u)}(p)$. Note that $L$ has one element in the tree of every node at scale $p$.

Case (i) $m \geq p$. Since every element of $L$ has proximity at most $p - 1$ with all other elements, $a_m = \sigma^p$ which is the maximum value it can take.

Case (ii) $m < p$ (assuming $p > 0$). Consider an arbitrary ordering of nodes at scale $m + 1$. We refer to the $q^{\text{th}}$ node in this ordering as "the $q^{\text{th}}$ node at scale $m + 1$".

Let the number of elements of $L$ belonging to the sub-tree of the $q^{\text{th}}$ node at scale $m + 1$ be $g_q, q = 0, \ldots, \sigma^{m+1} - 1$. We have

(7.38) $\quad a_m = \sum_{q=0}^{\sigma^{m+1}-1} g_q(\sigma^p - g_q) = \dfrac{\sigma^{2p+1+m}}{4} - \sum_{q=0}^{\sigma^{m+1}-1} (g_q - \sigma^p/2)^2$

since every element of $L$ in the tree of the $q^{\text{th}}$ node at scale $m + 1$ must have proximity *at most $m$* with all nodes *not* in the same tree but must have proximity *at least $m + 1$* with all nodes *within* the same tree.

The choice of $g_q$'s is constrained to lie on the hyperplane $\sum_q g_q = \sigma^p$. Obviously the quadratic form of (7.38) is maximized by the point on this hyperplane closest to the point $(\sigma^p/2, \ldots, \sigma^p/2)$ which is $(\sigma^{p-m-1}, \ldots, \sigma^{p-m-1})$. This is clearly achieved by $L \in \Gamma^{(u)}(p)$.

Now we prove the results for $L \in \Gamma^{(c)}(p)$.

Case (i) $m < D - p$. We have $a_m = 0$, the smallest value it can take.

Case (ii) $D - p \leq m < D$. Consider leaf node $\ell_i \in L$ which without any loss of generality belongs to the tree of first node at scale $m+1$. Let $a_m(\ell_i)$ be the number of elements of $L$ to which $\ell_i$ has proximity less than or equal to $m$. Now since $\ell_i$ has proximity less than or equal to $m$ only with those elements of $L$ not in the same tree, we must have $a_m(\ell_i) \geq \sigma^p - \sigma^{D-m-1}$. Since $L \in \Gamma^{(c)}(p)$ achieves this lower bound for $a_m(\ell_i), \forall i$ and $a_m = \sum_i a_m(\ell_i)$, $L \in \Gamma^{(c)}$ minimizes $a_m$ in turn. $\square$

We now study to what extent the above results transfer to the actual matrix of interest $Q_L^{-1}$. We start with a useful lemma.

**Lemma 7.4.** *Denote the eigenvalues of $Q_L$ by $\lambda_j, j = 1, \ldots, \sigma^p$. Assume that no leaf node of the tree can be expressed as a linear combination of other leaf nodes, implying that $\lambda_j > 0, \forall j$. Set $\mathcal{D}_L = [d_{i,j}]_{\sigma^p \times \sigma^p} := Q_L^{-1}$. Then there exist positive numbers $f_i$ with $f_1 + \ldots + f_p = 1$ such that*

$$\text{(7.39)} \qquad \sum_{i,j=1}^{\sigma^p} q_{i,j} = \sigma^p \sum_{j=1}^{\sigma^p} f_j \lambda_j, \text{ and}$$

$$\text{(7.40)} \qquad \sum_{i,j=1}^{\sigma^p} d_{i,j} = \sigma^p \sum_{j=1}^{\sigma^p} f_j / \lambda_j.$$

*Furthermore, for both special cases, $L \in \Gamma^{(u)}(p)$ and $L \in \Gamma^{(c)}(p)$, we may choose the weights $f_j$ such that only one is non-zero.*

*Proof.* Since the matrix $Q_L$ is real and symmetric there exists an orthonormal eigenvector matrix $U = [u_{i,j}]$ that diagonalizes $Q_L$, that is $Q_L = U \Xi U^T$ where $\Xi$ is diagonal with eigenvalues $\lambda_j, j = 1, \ldots, \sigma^p$. Define $w_j := \sum_i u_{i,j}$. Then

$$\text{(7.41)} \qquad \begin{aligned} \sum_{i,j} q_{i,j} &= \mathbf{1}_{1 \times \sigma^p} Q_L \mathbf{1}_{\sigma^p \times 1} = (\mathbf{1}_{1 \times \sigma^p} U) \Xi (\mathbf{1}_{1 \times \sigma^p} U)^T \\ &= [w_1 \ldots w_{\sigma^p}] \Xi [w_1 \ldots w_{\sigma^p}]^T = \sum_j \lambda_j w_j^2. \end{aligned}$$

Further, since $U^T = U^{-1}$ we have

$$\text{(7.42)} \qquad \sum_j w_j^2 = (\mathbf{1}_{1 \times \sigma^p} U)(U^T \mathbf{1}_{\sigma^p \times 1}) = \mathbf{1}_{1 \times \sigma^p} I \mathbf{1}_{\sigma^p \times 1} = \sigma^p.$$

Setting $f_i = w_i^2 / \sigma^p$ establishes (7.39). Using the decomposition

$$\text{(7.43)} \qquad Q_L^{-1} = (U^T)^{-1} \Xi^{-1} U^{-1} = U \Xi^{-1} U^T$$

similarly gives (7.40).

Consider the case $L \in \Gamma^{(u)}(p)$. Since $L = [\ell_i]$ consists of a symmetrical set of leaf nodes (the set of proximities between any element $\ell_i$ and the rest does not depend on $i$) the sum of the covariances of a leaf node $\ell_i$ with its fellow leaf nodes does not depend on $i$, and we can set

$$(7.44) \qquad \lambda^{(u)} := \sum_{j=1}^{\sigma^p} q_{i,j}(L) = c_D + \sum_{m=1}^{p} \sigma^{p-m} c_m.$$

With the sum of the elements of any row of $Q_L$ being identical, the vector $\mathbf{1}_{\sigma^p \times 1}$ is an eigenvector of $Q_L$ with eigenvalue $\lambda^{(u)}$ equal to (7.44).

Recall that we can always choose a basis of orthogonal eigenvectors that includes $\mathbf{1}_{\sigma^p \times 1}$ as the first basis vector. It is well known that the rows of the corresponding basis transformation matrix $U$ will then be exactly these normalized eigenvectors. Since they are orthogonal to $\mathbf{1}_{\sigma^p \times 1}$, the sum of their coordinates $w_j$ ($j = 2, \ldots, \sigma^p$) must be zero. Thus, all $f_i$ but $f_1$ vanish. (The last claim follows also from the observation that the sum of coordinates of the normalized $\mathbf{1}_{\sigma^p \times 1}$ equals $w_1 = \sigma^p \sigma^{-p/2} = \sigma^{p/2}$; due to (7.42) $w_j = 0$ for all other $j$.)

Consider the case $L \in \Gamma^{(u)}(p)$. The reasoning is similar to the above, and we can define

$$(7.45) \qquad \lambda^{(c)} := \sum_{j=1}^{\sigma^p} q_{i,j}(L) = c_D + \sum_{m=1}^{p} \sigma^m c_{D-m}.$$

$\square$

*Proof of Theorem 4.2.* Due to the special form of the covariance vector $\text{cov}(L, V_\emptyset) = \rho \mathbf{1}_{1 \times \sigma^k}$ we observe from (4.2) that minimizing the LMMSE $\mathcal{E}(V_\emptyset | L)$ over $L \in \Lambda_\emptyset(n)$ is equivalent to maximizing $\sum_{i,j} d_{i,j}(L)$ the sum of the elements of $Q_L^{-1}$.

Note that the weights $f_i$ and the eigenvalues $\lambda_i$ of Lemma 7.4 depend on the arrangement of the leaf nodes $L$. To avoid confusion, denote by $\lambda_i$ the eigenvalues of $Q_L$ for an arbitrary fixed set of leaf nodes $L$, and by $\lambda^{(u)}$ and $\lambda^{(c)}$ the only relevant eigenvalues of $L \in \Gamma^{(u)}(p)$ and $L \in \Gamma^{(c)}(p)$ according to (7.44) and (7.45).

Assume a positive correlation progression, and let $L$ be an arbitrary set of $\sigma^p$ leaf nodes. Lemma 7.3 and Lemma 7.4 then imply that

$$(7.46) \qquad \lambda^{(u)} \leq \sum_j \lambda_j f_j \leq \lambda^{(c)}.$$

Since $Q_L$ is positive definite, we must have $\lambda_j > 0$. We may then interpret the middle expression as an expectation of the positive "random variable" $\lambda$ with discrete law given by $f_i$. By Jensen's inequality,

$$(7.47) \qquad \sum_j (1/\lambda_j) f_j \geq \frac{1}{\sum_j \lambda_j f_j} \geq \frac{1}{\lambda^{(c)}}.$$

Thus, $\sum_{i,j} d_{i,j}$ is minimized by $L \in \Gamma^{(c)}(p)$; that is, clustering of nodes gives the worst LMMSE. A similar argument holds for the negative correlation progression case which proves the Theorem. $\square$